\documentstyle{article}
\input{amssymb.sty}

\date{}
\begin{document}
\title {Spectra of anticommutator  for two orthogonal projections
\thanks{ Corresponding author}
\thanks{ This
research was partially supported by the National Natural Science
Foundation of China(No.11571211, 11471200).}
\thanks{$E-mail\   address:$  hkdu@snnu.edu.cn(H.K.Du)} }
\author{Yan-Ni Dou$^{a}$, Yue-Qing Wang$^{b}$, Miao-Miao Cui$^{a}$ and Hong-Ke Du$^*$$^{a}$ }
 \maketitle\begin{center}
\begin{minipage}{12cm}
{\small a. College \ of \ Mathematics \ and \ Information \ Science,\
Shaanxi \ Normal \ University, Xi'an, 710062, P. R. China.  \\ b.Department of Mathematics and Physics, Chongqing University of Science and Technology, Chongqing, 401331, P. R. China. }
\end{minipage}
\end{center}
\begin{center}
\begin{minipage}{120mm}

{\bf Abstract}\\

{In this note, for any two orthogonal projections $P,Q$ on a Hilbert space, the characterization of spectrum of anticommutator $PQ+QP$ has been obtained. As a corollary, an alternative proof of the norm formula $$\parallel PQ+QP\parallel=\parallel PQ\parallel+\parallel PQ\parallel^2$$ has been got (see Sam Waltrs, Anticommutator norm formula for projection operators, arXiv:1604.00699vl [math.FA] 3 Apr 2016).
 }\\

\thanks {$AMS \  classification:$ \ 47A05}\\

\thanks{$Keywords \  and \  phrases:$   Orthogonal projection,
Anticommutator, Norm, Spectrum}
\end{minipage}
\end{center}

\vskip 1in \baselineskip 0.1in

\section{Introduction}

Throughout this note, a complex
Hilbert space is denoted by $\mathcal {H}$ , the algebra of all bounded linear
operators on $\mathcal {H}$ is denoted by $\cal B(H).$ For an operator $T\in {\cal B(H)}$,
$T^*$, $\sigma(T),$ $r_\sigma(T),$ ${\mathcal {N}}(T)$ and ${\mathcal {R}}(T)$
denote the adjoint, the spectrum, the spectral radius, the null space and the range of
$T$, respectively. $P\in \mathcal{B(H)}$ is said to be an orthogonal projection if $P=P^*=P^2.$ $A\in \mathcal{B(H)}$ is said to be self-adjoint if $A=A^*.$ $A\in \mathcal{B(H)}$ is said to be positive if $(Ax,x)\geq 0$ for each $x\in \mathcal{H}.$ $A\in \mathcal{B(H)}$ is said to be a contraction if $\parallel A\parallel\leq 1.$ The set of all contractions is denoted by ${\mathcal{B(H)}}_1.$

For two orthogonal projections $P$ and $Q,$ $PQ+QP$ is the anticommutator of $P$ and $Q,$ the main result of this note is proving the identity involving $\sigma(PQ+QP)$ and $\sigma(PQP).$

{\bf Theorem 1.1.} If $P$ and $Q$ are orthogonal projections on $\mathcal{H},$ then
\begin{equation}\small \{\lambda\pm \lambda^\frac{1}{2}:\lambda \in \sigma(PQP)\setminus \{0, 1\}\}\subseteq\sigma(PQ+QP)\subseteq\{\lambda\pm \lambda^\frac{1}{2}:\lambda \in \sigma(PQP)\}\bigcup \{0, 2\}.\end{equation}

As a consequence, we give an alternative proof of Theorem 1.1 in [3].

{\bf Corollary 1.2.} (see Theorem 1.1 in [3] or Theorem 1.3 in [4]) If $P$ and $Q$ are two orthogonal projections on Hilbert space $\mathcal{H},$ then $$\parallel PQ+QP\parallel=\parallel
PQ\parallel^2+\parallel PQ\parallel.$$

{\bf Proof.} The proof shall be divided into two cases.

First, if $ {\mathcal{R}}(P)\cap {\mathcal{R}}(Q)\neq \{0\},$ then $\parallel PQ+QP\parallel=2,$ it follows that $\parallel PQ\parallel=1$ and $\parallel QP\parallel=1.$ So $\parallel
PQ\parallel^2+\parallel PQ\parallel=2.$ Hence $\parallel PQ+QP\parallel=\parallel
PQ\parallel^2+\parallel PQ\parallel.$

Second, if  $ {\mathcal{R}}(P)\cap {\mathcal{R}}(Q)= \{0\},$
since $PQ+QP$ is self-adjoint and from the proof of Theorem 1.1 as below, then we have $$\begin{array}{rcl}\parallel PQ+QP\parallel&=&r_\sigma(PQ+QP)\\&=&r_\sigma(W) \\&=&
\max\{\mid\lambda\pm\lambda^\frac{1}{2}\mid: \lambda\in \sigma(PQP)\}\\&=&\max\{\lambda+\lambda^\frac{1}{2}: \lambda\in \sigma(PQP)\}\\&=&\max\{\lambda: \lambda\in \sigma(PQP)\}+\max\{\lambda^\frac{1}{2}: \lambda\in \sigma(PQP)\}\\&=&\parallel PQP\parallel+\parallel PQP\parallel^\frac{1}{2}\\&=&\parallel PQ\parallel^2+\parallel PQ\parallel.\end{array}$$

\section{ Proof of Theorem 1.1}

In this section, the proof of Theorem 1.1 has been based on Halmos modula of two projections theory. For the details of two projections theory, we refer the reader to [1]

To complete the proof of Theorem 1.1, as a preliminary, we begin with some lemmas.

 {\bf Lemma 2.1.} (see [1], [2] and [5])
If $P$ and $Q$ are two orthogonal projections on
$\mathcal{H}$, then $P$ and $Q$ have the
operator matrices
\begin{equation}P=I_1\oplus I_2\oplus I_3\oplus 0I_4 \oplus 0I_5\oplus
0I_6\end{equation} and
\begin{equation}Q=I_1\oplus 0I_2\oplus
\left(\begin{array}{cc}
Q_0&Q_0^\frac{1}{2}(I_3-Q_0)^\frac{1}{2}D\\
D^*Q_0^\frac{1}{2}(I_3-Q_0)^\frac{1}{2}&D^*(I_3-Q_0)D
\end{array}\right)\oplus 0I_5\oplus I_6\end{equation}
with respect to the space decomposition
${\mathcal{H}}=\oplus_{i=1}^6{\mathcal{H}}_i$, respectively, where
${\mathcal{H}}_1={\mathcal{R}}(P)\cap {\mathcal{R}}(Q)$,
${\mathcal{H}}_2={\mathcal{R}}(P)\cap {\mathcal{N}}(Q)$,
${\mathcal{H}}_3={\mathcal{R}}(P)\ominus ({\mathcal{H}}_1\oplus {\mathcal{H}}_2),$
${\mathcal{H}}_4={\mathcal{N}}(P)\ominus ({\mathcal{H}}_5\oplus {\mathcal{H}}_6),$
${\mathcal{H}}_5={\mathcal{N}}(P)\cap{\mathcal{N}}(Q)$ and ${\mathcal{H}}_6={\mathcal{N}} (P)\cap
{\mathcal{R}}(Q),$ $Q_0$ is a positive contraction on
${\mathcal{H}}_3$, $0$ and $1$ are not eigenvalues of $Q_0$, and $D$
is a unitary from ${\mathcal{H}}_4$ onto ${\mathcal{H}}_3$. $I_i$ are
the identity on ${\mathcal{H}}_i$, $i=1,\ldots,6.$

{\bf Lemma 2.2.} If $A\in \mathcal{B(H)}$ is a self-adjoint operator, then $\parallel A\parallel=r_\sigma(A).$

{\bf Proof of Theorem 1.1.}
If $P$ and $Q$ have operator matrices (2) and (3), respectively, then  \begin{equation}\small\begin{array}{rcl}&&PQ+QP\\&=&2I_1\oplus 0I_2\oplus
\left(\begin{array}{cc}
2Q_0&Q_0^\frac{1}{2}(I_3-Q_0)^\frac{1}{2}D\\
D^*Q_0^\frac{1}{2}(I_3-Q_0)^\frac{1}{2}&0
\end{array}\right)\oplus 0I_5\oplus 0I_6.\end{array}\end{equation}

Denote $$W=\left(\begin{array}{cc}
2Q_0&Q_0^\frac{1}{2}(I_3-Q_0)^\frac{1}{2}D\\
D^*Q_0^\frac{1}{2}(I_3-Q_0)^\frac{1}{2}&0
\end{array}\right).$$ We shall show that $$\sigma(W)=\{\lambda\pm \lambda^\frac{1}{2}:\lambda\in \sigma(Q_0)\}.$$

 For convenience, the proof will be divided into two cases.

{\bf Case one.} $0\in \sigma(W)$ if and only if at least one of $0$ and $1$ is in $\sigma(Q_0).$

If $0\in \sigma(W),$ then there exists a sequence $(x_n,y_n)^T_{1\leq n <\infty}$ of unit vectors with $x_n\in {\mathcal{H}}_3, y_n\in {\mathcal{H}}_4$ and $\parallel x_n\parallel^2+\parallel y_n\parallel^2=1, 1\leq n<\infty,$ such that $$\begin{array}{ll}&\left(\begin{array}{cc}
2Q_0&Q_0^\frac{1}{2}(I_3-Q_0)^\frac{1}{2}D\\
D^*Q_0^\frac{1}{2}(I_3-Q_0)^\frac{1}{2}&0
\end{array}\right)\left(\begin{array}{c}
x_n\\
y_n
\end{array}\right)\\=&\left(\begin{array}{c}
2Q_0x_n+Q_0^\frac{1}{2}(I_3-Q_0)^\frac{1}{2}Dy_n\\
D^*Q_0^\frac{1}{2}(I_3-Q_0)^\frac{1}{2}x_n
\end{array}\right)\rightarrow 0.\end{array} $$ If $0$ and $1$ are not in $\sigma(Q_0),$ then $Q_0$ and $I_3-Q_0$ are invertible. This implies that $x_n\rightarrow 0$ and $y_n\rightarrow 0.$ It is a contradiction.

Conversely, if at least one of $0$ and $1$ is in $\sigma(Q_0),$ then there exists a sequence $\{z_n\}\subset {\mathcal{H}}_4$ of unite vectors such that $$ Q_0^\frac{1}{2}(I_3-Q_0)^\frac{1}{2}Dz_n\rightarrow 0, \hbox{ as } n\rightarrow\infty,$$ this shows that $$\left(\begin{array}{cc}
2Q_0&Q_0^\frac{1}{2}(I_3-Q_0)^\frac{1}{2}D\\
D^*Q_0^\frac{1}{2}(I_3-Q_0)^\frac{1}{2}&0
\end{array}\right)\left(\begin{array}{c}0\\z_n\end{array}\right)
\rightarrow 0.$$ So $0\in \sigma(W).$

{\bf Case two.} $\sigma(W)\setminus \{0\}=(\sigma(Q_0+Q_0^\frac{1}{2})\cup\sigma(Q_0-Q_0^\frac{1}{2}))
\setminus \{0\}.$

If $\lambda\in \sigma(W)$ with $\lambda\neq 0,$ then, by Schur complement, $$W-\lambda=\left(\begin{array}{cc}
2Q_0-\lambda I_3&Q_0^\frac{1}{2}(I_3-Q_0)^\frac{1}{2}D\\
D^*Q_0^\frac{1}{2}(I_3-Q_0)^\frac{1}{2}&-\lambda I_4
\end{array}\right)$$ is not invertible if and only if so is $(\lambda-(Q_0+Q_0^\frac{1}{2}))
(\lambda-(Q_0-Q_0^\frac{1}{2}))$ since $$\begin{array}{rcl}&&2Q_0-\lambda+Q_0^\frac{1}{2}
(I_3-Q_0)^\frac{1}{2}D\lambda^{-1}
D^*
Q_0^\frac{1}{2}(I_3-Q_0)^\frac{1}{2}
\\&=&\lambda^{-1}(\lambda(2Q_0-\lambda)+Q_0(I_3-Q_0))\\&=&
-\lambda^{-1}(\lambda-(Q_0+Q_0^\frac{1}{2}))
(\lambda-(Q_0-Q_0^\frac{1}{2})).\end{array}$$ This shows that at leat one of $\lambda-(Q_0+Q_0^\frac{1}{2})$ and
$\lambda-(Q_0-Q_0^\frac{1}{2})$ is not invertible.
Hence,
$$\sigma(W)\setminus \{0\}=(\sigma(Q_0+Q_0^\frac{1}{2})\cup\sigma(Q_0-Q_0^\frac{1}{2}))
\setminus \{0\}.$$

Combining Case two with Case one, we get $$\sigma(W)=\sigma(Q_0+Q_0^\frac{1}{2})\cup\sigma(Q_0-Q_0^\frac{1}{2})
=\{\lambda\pm \lambda^\frac{1}{2}:\lambda\in \sigma(Q_0)\}.$$

 Next, we shall consider the structure of $\sigma(PQP)$ and $\sigma(PQ+QP).$ On the sake of convenience, the rest of the proof has been divided in two cases.

 Case one. Assume that ${\mathcal{H}}_1\neq \{0\}.$ In this time,
since $PQP=I_1\oplus 0I_2\oplus Q_0\oplus 0I_4\oplus 0I_5\oplus 0I_6$ and $PQ+QP=2I_1\oplus 0I_2\oplus W\oplus 0I_5\oplus 0I_6,$ we have $$\sigma(Q_0)\subseteq \sigma(PQP)\subseteq \sigma(Q_0)\bigcup \{0, 1\}$$ and $$\sigma(W)\subseteq \sigma(PQ+QP)\subseteq \sigma(W)\bigcup\{2,0\}.$$ So that, $$\small \{\lambda\pm \lambda^\frac{1}{2}:\lambda \in \sigma(PQP)\setminus \{0, 1\}\}\subseteq\sigma(PQ+QP)\subseteq\{\lambda\pm \lambda^\frac{1}{2}:\lambda \in \sigma(PQP)\}\bigcup \{0, 2\}.$$

Case two. Assume that ${\mathcal{H}}_1= \{0\}.$ We have $$\sigma(Q_0)\subseteq \sigma(PQP)\subseteq \sigma(Q_0)\bigcup \{0\}$$ and $$\sigma(W)\subseteq \sigma(PQ+QP)\subseteq \sigma(W)\bigcup\{0\}.$$
Hence, $$\small \{\lambda\pm \lambda^\frac{1}{2}:\lambda \in \sigma(PQP)\setminus \{0\}\}\subseteq\sigma(PQ+QP)\subseteq\{\lambda\pm \lambda^\frac{1}{2}:\lambda \in \sigma(PQP)\}\bigcup \{0\}.$$

Finally, combining Case two and Case one, we obtain $$\small \{\lambda\pm \lambda^\frac{1}{2}:\lambda \in \sigma(PQP)\setminus \{0, 1\}\}\subseteq\sigma(PQ+QP)\subseteq\{\lambda\pm \lambda^\frac{1}{2}:\lambda \in \sigma(PQP)\}\bigcup \{0, 2\}.$$

{\bf Remark.} Since $0\in \sigma(W)$ if and only if at least one of $0\in \sigma(Q_0)$ and $1\in \sigma(Q_0)$ holds, so $$\sigma(W)=\sigma(Q_0+Q_0^\frac{1}{2})\bigcup
\sigma(Q_0-Q_0^\frac{1}{2}).$$

{\bf Corollary 2.3.} If $P$ and $Q$ are two orthogonal projections, then $$\sigma(PQ+QP)\subseteq [-\frac{1}{4}, 2].$$

{\bf Proof.} Observing that $\sigma(Q_0)\subseteq [0, 1]$ and $$\begin{array}{rcl}\{\lambda\pm \lambda^\frac{1}{2}:\lambda\in \sigma(Q_0)\}&\subseteq &\sigma(PQ+QP)\\&\subseteq &\{\lambda\pm \lambda^\frac{1}{2}:\lambda \in \sigma(PQP)\bigcup \{0, 1\}\}\\&\subseteq &\{\lambda\pm \lambda^\frac{1}{2}:\lambda \in [0, 1]\},\end{array}$$
 we get $$\sigma(PQ+QP)\subseteq [-\frac{1}{4}, 2].$$

{\bf Example 2.4.} There exists a pair $(P,Q)$ of orthogonal projections with $$\sigma(PQ+QP)=[-\frac{1}{4}, 2].$$

Let $\{e_i\}_{1\leq i<\infty}$ be an orthonormal basis of Hilbert space $\mathcal{K}$ and let $\{r_i\}_{1\leq i<\infty}$ the set of all rational numbers in $(0, 1).$ Define a positive contraction $Q_0$ on $\mathcal{K}$ by $$Q_0e_i=r_ie_i, 1\leq i<\infty,$$ it is obvious that $\sigma(Q_0)=[0,1],$ and define $P$ and $Q$ on ${\mathcal{H}}={\mathcal{K}}\oplus{\mathcal{K}}$ by $$P=\left(\begin{array}{cc}I_{\mathcal{K}}&0\\0&0\end{array}\right) \hbox{ and } Q=\left(\begin{array}{cc}Q_0&Q_0^\frac{1}{2}
(I_{\mathcal{K}}-Q_0)^\frac{1}{2}\\Q_0^\frac{1}{2}
(I_{\mathcal{K}}-Q_0)^\frac{1}{2}&I_{\mathcal{K}}-Q_0\end{array}\right).$$
It is clear that $P$ and $Q$ are orthogonal projections and $$\sigma(PQ+QP)=[-\frac{1}{4}, 2].$$

In fact, if $\triangle$ is a closed subset of $[-\frac{1}{4}, 2],$ then there exists a pair $(P, Q)$ of orthogonal projections such that $\sigma(PQ+QP)=\triangle.$

\end{document}